\documentclass[12pt,thmsa,emstex]{article}
\UseRawInputEncoding
\usepackage{amsfonts}
\usepackage{t1enc}
\usepackage{amssymb}
\usepackage{epsfig}
\usepackage[french,english]{babel}
\usepackage{amsmath}
\usepackage{graphics}
\usepackage{layout}
\usepackage{latexsym}
\usepackage{color}
\usepackage{verbatim}

\usepackage{xcolor} 
\long\def\metanote#1#2{{\color{#1}\ 
\ifmmode\hbox\fi{\sffamily\mdseries\upshape [#2]}\ }} 

\def\mylabel#1{\label{#1}}

\newtheorem{theorem}{Theorem}[section]
\newtheorem{lemma}[theorem]{Lemma}
\newtheorem{corollary}[theorem]{Corollary}
\newtheorem{proposition}[theorem]{Proposition}

\newtheorem{exercise}[theorem]{Exercise}

\newtheorem{remark}[theorem]{Remark}
\newtheorem{example}[theorem]{\bf{Example}}
\newtheorem{assumption}[theorem]{\bf{Assumption}}

\def\bit{\begin{itemize}}
\def\eit{\end{itemize}}

\reversemarginpar   

\def\bc{\begin{center}}
\def\ec{\end{center}}

\def\bthm{\begin{theorem}}
\def\ethm{\end{theorem}}

\def\bcor{\begin{corollary}}
\def\ecor{\end{corollary}}

\def\bprop{\begin{proposition}}
\def\eprop{\end{proposition}}

\def\blem{\begin{lemma}}
\def\elem{\end{lemma}}

\def\bex{\begin{example}}
\def\eex{\end{example}}

\def\bexo{\begin{exercise}}
\def\eexo{\end{exercise} }

\def\brem{\begin{remark}}
\def\erem{\end{remark}}
\def\prf{{\bf Proof: }}

\def\bdes{\begin{description}}
\def\edes{\end{description}}

\def\iti{\item[(i)]}
\def\itii{\item[(ii)]}

\def\itiii{\item[(iii)]}

\def\beq{\begin{equation}}
\def\eeq{\end{equation}}

\def\ben{\begin{enumerate}}
\def\een{\end{enumerate}}

\def\beqar{\begin{eqnarray}}
\def\eeqar{\end{eqnarray}}

\def\beqarr{\begin{eqnarray*}}
\def\eeqarr{\end{eqnarray*}}

\def\qed{\hfill $\Box$ \\[2ex]}
\def\prf{{\bf Proof: }\hspace{.1in}}

\catcode`\@=11

\newcommand{\M}{\mathcal{M}}

\newcommand{\EE}{\mathbb{E}}
\newcommand{\E}{\mathsf{E}}

\def\Ind{{\mathbf 1}}
\def\RR{{\mathbb R}}  
\def\Rp{{\mathbb R}_+}   

\def\NN{{\mathbb N}}

\def\Pr{{\mathsf P}}
\def\M{{\mathbf M}}
\def\m{{\mathbf m}}
\def\i{{\mathbf i}}
\def\j{{\mathbf j}}

\def\rar{\rightarrow}
\def\eps{\varepsilon}

\begin{document}
\title{Regularity of invariant densities for random switching between two linear odes in $\RR^d.$}
\author{Michel  Bena{\"i}m$^1$  and Am\'ethyste Bichard$^2$}
\footnotetext[1]{Institut de Math\'{e}matiques, Universit\'{e} de Neuch\^{a}tel, Switzerland.}
\footnotetext[2]{ENS, Paris Saclay.}
\maketitle

\begin{abstract}
In a paper entitled {\em singularities of invariant densities for random switching between two linear odes in $2D$}, Bakhtin et al \cite{BHLM21}, consider a Markov process obtained by random switching between two  stable linear vector fields in the plane  and  characterize the singularities of the invariant density in terms of the switching and contraction rates. This paper considers a generalization of this  model obtained by random switching between two  stable linear vector fields in $\RR^d$ and provides sufficient conditions ensuring that the invariant distribution is absolutely continuous and  has a $C^r$ density. In dimension $\geq 3$ it provides,  to the best of our knowledge, the first fully non-elliptic  example of random switching for which quantitative conditions guaranteeing smoothness of the invariant density can be proved.
\end{abstract}
\tableofcontents
\section{Introduction}
Processes obtained by random switching between finitely many ODEs provide an interesting class of {\em piecewise deterministic Markov processes} (PDMPs), a term coined by Davis \cite{Davis84}, that are often used in a variety of fields including molecular and cellular biology \cite{Bressloff17}, \cite{Lawley15}, population dynamics, ecology \cite{BenaimLobry}, \cite{MR3695481}, \cite{MR3582809},  MCMC simulation \cite{MR3983339}, \cite{Monemvassitis_2023}, and elsewhere \cite{Yin}.
The investigation of their ergodic properties started around a decade ago (\cite{bakhtin&hurt}, \cite{BLMZ_2015},  \cite{BLMZ_2019}, \cite{Cloez15},  \cite{BCL17}, \cite{BHS18}). One of the key results established in \cite{bakhtin&hurt} and \cite{BLMZ_2015} (see also \cite{BenaimHurth} chapter 6.4, and \cite{benaïm2024invariant}, Section 4.2) is that, under the existence of an {\em accessible} and {\em hypoelliptic} point,  there is  (at most) one invariant distribution which is absolutely continuous with respect to the Lebesgue measure, and whose support equates the set of accessible points.  Here, by hypoelliptic, we mean a point at which the Lie algebra generated by the vector fields has full rank. By accessible, we mean a point whose every neighborhood can be reached from everywhere by following the  solution curves of the different ODEs.  It is interesting to point out that this type of result has long been known for stochastic differential equations. This was proved, for example, by Arnold and Kiemann in \cite{ArnoldKliemann}.  An alternative proof, relying on classical results by Bony \cite{Bony} is in  \cite{BenaimHurth}, Theorem 6.34. However, while for stochastic differential equations, hypoellipticity ensures smoothness of the invariant densities (by Hormand\"er's theorem), the situation is more involved for PDMPs obtained by random switching. The smoothness effect induced by hypoellipticity may be balanced by contracting properties of the ODEs.
In dimension one, thanks to the work of Bakhtin, Hurth and  Mattingly \cite{bakhtin&hurt&matt}, the situation is well understood:   Singularities may occur at attracting equilibria of the vector fields, while the invariant densities are always smooth away from these equilibria.   In higher dimensions, the situation is much more complicated and, with the exception of a few very specific examples, little was known until recently.
In \cite{BHLM18}, Bakhtin, Hurth, Lawley and  Mattingly investigated a  two dimensional {\em elliptic} (meaning that the vector fields are transverse)  system on the torus given by two  free periodic vector fields . They showed that, under certain strong assumptions, the densities are smooth. In \cite{BHLM21}, the same authors considered a planar system given by two stable linear vector fields, which is elliptic  almost everywhere, except on a line where it is "hypoelliptic".  They were able to provide precise conditions (on the jump rates and the eigenvalues  of the vector fields) ensuring boundedness  and unboundedness of the densities. This example is one of the main motivations of the present paper.
The two recent articles \cite{BT23} and \cite{benaïm2024invariant} consider the more general question of proving / disproving the $C^r$ regularity of invariant densities in any dimension.  The main result of \cite{BT23} (see also Section 4.4 in \cite{benaïm2024invariant}) is that under a strong version of hypoellipticity called the {\em $1$-Hörmander condition} (meaning that the family of the vector fields and their first-order brackets has full rank), densities are $C^k$ provided the jump rate is fast enough. The paper \cite{benaïm2024invariant} provides some quantitative estimates and explains how certain dynamical quantities (like the expansion rates and the  expansion volume rate) should be taken into account to ensure $C^k$ regularity. This was successfully applied to the elliptic system on the torus considered in \cite{BHLM18} under much weaker hypotheses.

The present paper  builds on the ideas introduced in \cite{benaïm2024invariant} and gives precise conditions ensuring  that the invariant probability measure of a system obtained by random switching between two linear ODEs in $\RR^d$ possesses a $C^r$ density. This model is a natural generalization of the planar example considered in \cite{BHLM21}. In dimension $\geq 3$ it provides,  to the best of our knowledge, the first fully non-elliptic  example of random switching for which quantitative condition guaranteeing smoothness of the invariant density can be proved. Note also that, for this example,  the stronger $1$-Hormandër condition used in \cite{BT23} is not verified.
\subsection{Model and examples}
Throughout $d$ is an integer  $\geq 2.$ The Euclidean scalar product on $\RR^d$ is denoted $\langle \cdot \rangle$ and the associated norm  $\|\cdot\|.$ Let $A$ be a $d \times d$ real matrix whose eigenvalues have negative real parts:
$$\Lambda_1 \leq \ldots \leq \Lambda_d < 0.$$
Let $p \in \RR^d \setminus \{0\}.$ For $x \in \RR^d,$ set $F_0(x) = Ax$ and $F_1(x) = A (x-p).$ We consider the PDMP $(X_t,I_t)_{t  \geq 0}$ living on $\RR^d \times \{0,1\}$ defined as
\beq
\label{eq:pdmp}
\frac{dX_t}{dt} = F_{I_t}(X_t),
\eeq
where $(I_t)_{t \geq 0}$ is a Markov chain on $\{0,1\},$ independent of $(X_t)_{t \geq 0},$ with transition rates $\alpha_{01}, \alpha_{10} > 0.$ That is
\beq
\label{eq:jumps}
\Pr(I_{t+s} =  1-i| I_t = i, {\cal F}_t) = s \alpha_{i(1-i)}  + o(s)
\eeq
for $i \in \{0,1\},$ where ${\cal F}_t = \sigma\{(X_s,I_s); s  \leq t\}.$

Let $$H := \mathsf{span}(\{ A^k p: \: k \in \NN\}),$$
 be the vector space spanned by the family  $\{ A^k p: \: k \in \NN\}.$
It is easy to verify that
$H$ has dimension $d' \in \{1, \ldots, d\}$ where $d'$ is the largest integer such that $p, Ap, \ldots, A^{d'-1} p$ are linearly independent.

The next elementary lemma shows that  $H \times \{0,1\}$ is an  {\em attracting invariant} set for the PDMP ((\ref{eq:pdmp}), (\ref{eq:jumps})).
 For $x \in \RR^d,$ we let $P_H(x)$ denote the orthogonal projection of $x$ onto $H$ and $\mathsf{dist}(x,H) = \|x-P_H(x)\|.$
\blem
\label{lem:attract}
Let $0 < \lambda < |\Lambda_d|.$ Let  $(X_t,I_t)$ be a solution to ((\ref{eq:pdmp}), (\ref{eq:jumps})). Then
\bdes
\iti For all $t \geq 0,$
  $$X_t \in H \Leftrightarrow \{X_s \: : s \geq 0\} \subset H$$
  \itii  $$\mathsf{dist}(X_t,H) \leq C e^{-\lambda t} \mathsf{dist}(X_0,H)$$ where $C$ is some constant depending on $A.$
  \edes
  \elem
  \prf $(i)$ follows from the fact that both $F_0$ and $F_1$ are tangent to $H.$ For $(ii),$ let $(\tilde{X_t})_{t \geq 0}$ be the solution to (\ref{eq:pdmp}) starting from $\tilde{X_0} = P_H(X_0).$ Then $\tilde{X_t} \in H$ and
  $$\|X_t - \tilde{X}_t\| = \|e^{tA} (X_0 -\tilde{X}_0)\| \leq C e^{-\lambda t}\|(X_0 -\tilde{X}_0)\|.$$ \qed
 In view of this lemma, the long term behavior of the PDMP ((\ref{eq:pdmp}),(\ref{eq:jumps})) on $\RR^d \times \{0,1\}$ is given by its long term behavior on $H \times \{0,1\}.$ Therefore, by replacing $\RR^d$ by $H$ and $A$ by $A|_H$ if necessary, we can assume without loss of generality that the following condition is met.
\begin{assumption}[Standing assumption]
\label{mainhyp}
$H = \RR^d,$ or equivalently: $p, Ap, \ldots, A^{d-1} p$ are linearly independent.
\end{assumption}
From now on, this assumption is implicitly assumed to be verified.
 \bex
 \label{ex1}
 {\rm This example is given by Malrieu in   \cite{Malrieu}, Section 4.2.
 Here $d = 2,$
$$A =  \left(
   \begin{array}{cc}
     -1 & -1 \\
     1 & -1 \\
   \end{array}
 \right), \mbox{ and } p = \left(
                \begin{array}{c}
                  1 \\
                  0 \\
                \end{array}
              \right)
 .$$
Malrieu (\cite{Malrieu}, Section 4.2, Open question 4) asks: {\em  What can be said on the smoothness of the density of such a process ?}}
 \eex
 \bex
 \label{ex2}
 {\rm This example is the one thoroughly studied  by Bakhtin, Hurth, Lawley and Mattingly in \cite{BHLM21}. Here again $d = 2$ and
 $$A =  \left(
   \begin{array}{cc}
     -\lambda_1 & 0 \\
     0 & -\lambda_2 \\
   \end{array}
 \right), p = \left(
                \begin{array}{c}
                  1 \\
                  1 \\
                \end{array}
              \right),$$ with $\lambda_1 > \lambda_2 > 0.$
The results in \cite{BHLM21} characterize precisely the boundedness properties of the invariant distribution. It is shown that  this distribution is absolutely continuous with respect to Lebesgue (see Section \ref{sec:main} for a precise definition)  and that the  density is bounded on every compact set contained in the interior of the accessible set. At points $0$ and $p$ as well as at boundary points of the accessible set, this density  blows up  if the switching rates are sufficiently small, while it remains bounded if they are large enough. We refer the reader to Theorems 1 and 2 in \cite{BHLM21}.}
 \eex
\bex
\label{ex3} {\rm A natural generalization of Example \ref{ex2} in dimension $d \geq 2$ is given by $A = \mathsf{diag}(-\lambda_1, \ldots, -\lambda_d)$  with
$\lambda_1 > \ldots > \lambda_d > 0,$ and $p = (1, \ldots, 1)^t.$ Assumption \ref{mainhyp} is satisfied because the matrix $(p,Ap, \ldots, A^{d-1} p)$ is the Vandermonde matrix whose determinant is $\prod_{i < j} (\lambda_j-\lambda_i) \neq 0.$
}
\eex
\brem
\label{rem2d}
{\rm In dimension $2,$ (see e.g~ Examples \ref{ex1} and \ref{ex2}) the set of points at which the system is {\em nonelliptic}  - that is the set of  points $x \in \RR^2$ at which $(F_0(x), F_1(x))$ has rank $<2$ - is the one dimensional vector space $\RR p.$ On $\RR p$, the system satisfies the {\em $1$-bracket condition} :
$$\mathsf{span}(F_0(x), F_1(x), [F_0, F_1](x)) = \RR^2.$$
Therefore, Theorem 4.15 and Example 4.18 in \cite{benaïm2024invariant},  imply that for all $k \geq 0,$  the invariant density of the process is $C^k$ provided that $\alpha_{01}$ and $\alpha_{10}$ are chosen {\em sufficiently large}. Note, however, that  these are purely qualitative results that  give no information on the size of $\alpha_{01}$ and $\alpha_{10}.$
}
\erem
\brem
\label{rem3d}
{\rm
In dimension $d \geq 3$ , there are always points at which the $1$-bracket condition is not satisfied (in particular, the results mentioned in Remark \ref{rem2d}  relying on the $1$-bracket conditions are useless). Still, the system is hypoelliptic because
the family
$$([F_0,F_1], [F_0,[F_0,F_1]], \ldots [F_0, \ldots, [F_0,[F_0,F_1]]) = (Ap, A^2p, \ldots, A^d p)$$ has rank $d.$
}
\erem
\section{Main results}
\mylabel{sec:main}
Let $(P_t)_{t \geq 0}$ be the Markov semi-group associated to $(X_t,I_t)_{t \geq 0}.$
 That is $$P_t f(x,i) =  \E (f(X_t,I_t) | X_0 = x, I_0 = i),$$ for all $f : \RR^d \times \{0,1\} \rar \RR$ bounded and measurable.

 As usual, if $\mu$ is a measure on $\RR^d \times \{0,1\}$ and $f$ a measurable, bounded or nonnegative, function on $\RR^d \times \{0,1\},$ $\mu f$ stands for $\int f d \mu,$ and $\mu P_t$ is the measure defined by $$(\mu P_t) f = \mu (P_t f).$$ We also adopt the classical notation
 $P_t((x,i), \cdot)  = \delta_{x,i} P_t.$

 A probability measure $\Pi$ on $\RR^d \times \{0,1\}$ is said to be {\em invariant} for $(P_t)_{t \geq 0}$ (or $(X_t,I_t)_{t \geq 0}$) if $$\Pi P_t = \Pi$$ for all $t \geq 0.$

 The {\em Lebesgue measure} on $\RR^d \times \{0,1\}$ is the measure $\m = dx \otimes (\delta_0 + \delta_1)$ where $dx$ stands for the Lebesgue measure on $\RR^d.$

If a measure $\mu$ on $\RR^d \times \{0,1\}$ is absolutely continuous with respect to $\m$, its {\em density} is a $L^1(\m)$ function $\rho : \RR^d \times \{0,1\} \rar \Rp, (x,i) \mapsto \rho_i(x)$ where $\rho_i = \frac{d \mu_i}{dx},$ and $\mu_i(\cdot):=\mu (\cdot \times \{i\}).$ We say that $\rho$ is continuous, $C^k,$ etc. whenever $\rho_i, i = 0,1,$ are continuous, $C^k,$ etc.

 Given a  piecewise continuous functions  $j : \Rp \rar \{0,1\},$ and  $x \in \RR^d,$ let $t \mapsto y(t,x,j)$ denote the solution to the non-autonomous ODE
$$\frac{dy}{dt} = F_{j(t)}(y),$$ starting from $y(0) = x.$ Let $\gamma^+(x)$ be the set of all  $y(t,x,j)$ with $t \geq 0,$ and  $j$ piecewise continuous. The {\em accessible set} of the pair $F_0,F_1$ is the nonempty compact connected set (see Proposition 3.11 in \cite{BLMZ_2015} and \cite{BLMZ_2019} for more details) defined by
$$\Gamma = \bigcap_{x \in \RR^d} \overline{\gamma^+(x)}.$$
\begin{figure}
\centering
\includegraphics[width=12cm]{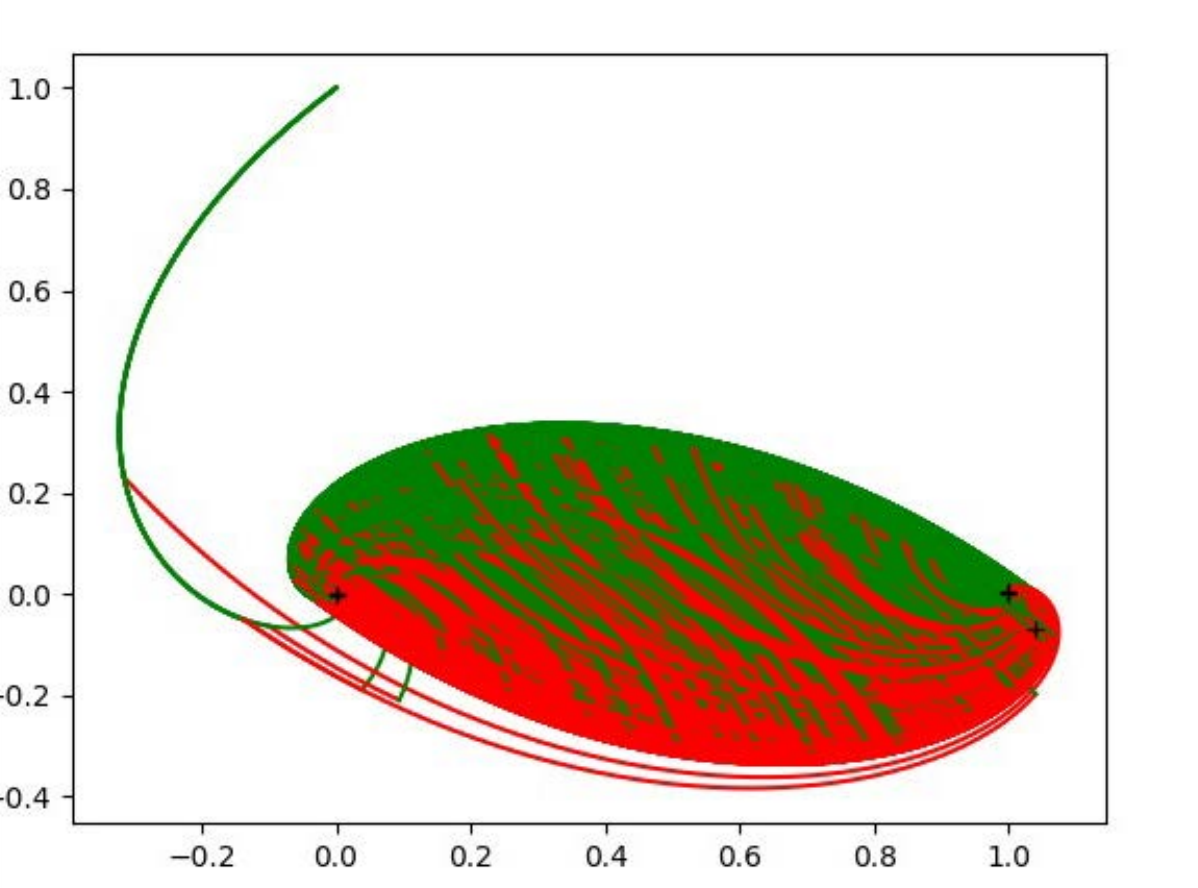}\\
\caption{The set $\Gamma$ in Example \ref{ex1} \label{Malrieu}}
\end{figure}
 The following proposition  easily follows from  existing results.
\bprop
\label{prop:background}
\bdes
\iti
The process $(X_t,I_t)_{t \geq 0}$ has a unique invariant probability measure $\Pi$ which is absolutely continuous with respect to $\m$ and whose density $\rho$ can be chosen to be lower semi-continuous.
\itii
The topological support of $\Pi$ writes $\Gamma \times \{0,1\}.$  Furthermore $$\mathsf{Int}(\Gamma) = \{x \in \RR^d \: : \rho_i(x) > 0\},$$
and $$\overline{\mathsf{Int}(\Gamma)} = \Gamma.$$
\itiii There exists $\gamma > 0$ and, for all $R > 0,$ a constant $C_R > 0$ such that for all $x \in B(0,R)$ (the Euclidean ball centered at the origin with radius $R$) and $i = 0,1,$
$$\|P_t((x,i), \cdot) - \Pi\|_{TV}  \leq  C_R e^{-\gamma t}.$$
Here $\|\cdot\|_{TV}$ stands for the total variation distance.
\edes
\eprop
\prf It is not hard to show (and this is done in Section \ref{sec:prelim}) that $(X_t)$ eventually enters  a compact set $M \subset \RR^d$ positively invariant for the flows induced both by $F_0$ and $F_1.$ Furthermore, as pointed out in Remark \ref{rem3d}, the Lie algebra generated by $F_0$ and $F_1$ has full rank. Statements $(i)$ and $(iii)$ can then be deduced from Corollary 2.7 in  \cite{BHS18} (itself relying on previous results in \cite{bakhtin&hurt}, \cite{BLMZ_2015} and \cite{li17}). The second statement follows from Theorem 4.5 in \cite{benaïm2024invariant}.
\qed
The main result of the paper is the following.
\bthm
\label{th:main}
Let $r \geq 0$ be a  integer. Assume that
\beq
\label{eq:maincond}
\min( \alpha_{01}, \alpha_{10} ) > -[\sum_{i = 1}^d \Lambda_i + r \Lambda_1].
 \eeq
 Then, $\rho$ is $C^r$ on $\RR^d \times \{0,1\},$ and $\{\rho > 0\} =  \mathsf{Int}(\Gamma),$ where $\rho$ is, like in Proposition \ref{prop:background}, the density of $\Pi.$
\ethm
\brem
{\rm The fact that $\{\rho > 0\} =  \mathsf{Int}(\Gamma)$ in Theorem \ref{th:main}, follows directly from Proposition \ref{prop:background}. Indeed, by this proposition, we know that there exists a lower semi-continuous  version of $\rho,$ say $\tilde{\rho},$ such that $\{\tilde{\rho} > 0\} =  \mathsf{Int}(\Gamma).$ The set $\{\tilde{\rho} - \rho > 0\}$ is empty, as an open set (by lower semi-continuity)  having zero measure. This shows that  $\rho > 0$ on  $\mathsf{Int}(\Gamma).$ The converse inclusion is obvious.
 }
 \erem
\brem
{\rm In case Assumption \ref{mainhyp} is not satisfied, the process $(X_t,I_t)_{t \geq 0}$ still admits a unique invariant probability measure $\Pi,$ supported by $H.$ Furthermore, it is not hard to show (by using a standard coupling argument), that for all $0 < \lambda < \min\{|\Lambda_d|, \alpha_{01} + \alpha_{10}\}$ there exists some constant $C > 0$ such that for every probability measure $\mu$ on $\RR^d \times \{0,1\}$
 $$W_1(\mu P_t, \Pi) \leq C e^{-\lambda t} W_1(\mu, \Pi)$$ where $W_1$ is the $1$-Wasserstein distance associated to the distance on $\RR^d \times \{0,1\}$ defined as $D((x,i),(y,j)) = \|x-y\| + \Ind_{i \neq j}.$

 Clearly the convergence of $(\mu P_t)$ to $\Pi$ cannot hold in total variation because, by invariance of $H,$ (see Lemma \ref{lem:attract}) $|\mu P_t(H) - \Pi(H)| = 1-\mu(H).$ However, for the process $(X_t,I_t)_{t  \geq 0}$ {\em restricted} to $H \times \{0,1\},$ Proposition \ref{prop:background} and Theorem \ref{th:main} hold true with $\RR^d$ replaced by $H$ and $\m$ the Lebesgue measure on $H \times \{0,1\}.$
}
\erem
\section{Proof of Theorem \ref{th:main}}
\label{sec:proof}
\subsection{Preliminaries}
\label{sec:prelim}
For $i = 0,1,$ we let $\{\Phi_i^t\}_{t \in \RR}$ denote the flow induced by $F_i$ on $\RR^d.$
That is,
$$\Phi_0^t(x) = e^{t A} x, \mbox{ and } \Phi_1^t(x) = e^{t A}(x-p) + p.$$
\subsubsection*{Adapted norms and state space}
The following lemma is folklore. We provide a proof for convenience.
\blem[Adapted norms]
\label{lem:norm}
For all  $\eps > 0,$ there exist  scalar products $( \,, )_{\eps}$ and $\langle\, , \rangle_{\eps}$  whose associated norms are denoted $N_{\eps}$  and $\|\:\|_{\eps}$ such that:
 $$( A x, x )_{\eps} \leq  (\Lambda_d + \eps) N_{\eps}^2(x)$$ and
 $$\langle A x, x \rangle_{\eps} \geq  (\Lambda_1 - \eps) \|x\|^2_{\eps}.$$
In particular, if $0 < \eps < |\Lambda_d|$, there exists $R > 0$  such that
$$\frac{d N_{\eps}(\Phi^t_i(x))}{dt} < 0$$ for all $i \in \{0,1\}, t \geq 0$ and $N_{\eps}(x) \geq  R.$
\elem
\prf  Recall that $\|\cdot\|$ and $\langle, \, \rangle$ stand for the Euclidean norm and scalar product. Fix $\eps > 0$ and let $\lambda = \Lambda_d + \eps.$ Then, there exists $\tau > 0,$ sufficiently large, so that $\|e^{tA} x\| \leq e^{\lambda t}\|x\|$ for all $t \geq \tau.$ Set  $$(x, y)_{\eps} = \int_0^{\tau} \langle e^{sA} x, e^{sA} y\rangle e^{-2 \lambda s}ds.$$ To shorten notation we write $(x,y)$ for  $(x,y)_{\eps},$ and $N$ for $N_{\eps}.$
Then $$(A x, x) =  \frac{1}{2} \int_0^{\tau} \frac{d \|e^{sA} x\|^2}{ds} e^{-2 \lambda s}ds = \frac{1}{2}(\|e^{\tau A} x\|^2 e^{-2\lambda\tau} -\|x\|^2) +\lambda N^2(x) \leq \lambda N^2(x).$$ This proves the first assertion.  The second is analogous. It immediately follows that whenever $\lambda < 0,$
 $(F_0(x),x) < 0$ for all $x \neq 0,$ and $$(F_1(x), x ) =  ( A(x-p), (x-p) ) + (A(x-p),p) \leq N(x-p)[ \lambda N(x-p) + C]$$ with $C = N(A) N(p).$ This latter quantity is negative for $N(x)$ sufficiently large. This proves the last statement.
\qed
We now define the {\em state space} $\M$ of the process. Fix  $0 < \eps_0 < |\Lambda_d|.$ The value of $\eps_0$ is unimportant. Let $N = N_{\eps_0}$ and  $\M = M \times \{0,1\}$ where
\beq
\label{eq:defM}
M = \{x \in \RR^d: \: \:N(x) \leq R\}
\eeq with $R$ given by the preceding lemma.
From what precedes, the PDMP $(X_t,I_t)$ starting from $(x,i) \in \RR^d \times \{0,1\}$ eventually enters $\M$ and remains in $\M$ afterwards.
\subsubsection*{The spaces $C^r(M), C^r(\M)$}
For $r \in \NN,$ and $M$ as above (\ref{eq:defM}), let $C^r(M)$ be the space of $C^r$ functions $f: \RR^d \rar \RR$ such  that $f = 0$ on $\RR^d \setminus M.$
For $\eps > 0,$ let $\|\cdot\|_{\eps}$ be the norm on $\RR^d$ induced by $\langle , \rangle_{\eps}$ as defined in Lemma \ref{lem:norm}. For $f \in C^r(M),$ define
\beq
\label{eq:normm} \|f\|_{r,\eps} = \sum_{k = 0}^r \|D^k f\|_{\eps},
\eeq
where
$D^0 f = f$ and  $\|D^k f\|_{\eps} = \sup_{x \in M} \|D^k f(x)\|_{\eps}.$ Here, as usual, for $k \geq 1,$
$$\|D^k f(x)\|_{\eps} = \sup \{|D^k f(x)(u_1, \ldots, u_k)|:  \: \|u_1\|_{\eps} \leq 1, \ldots, \|u_k\|_{\eps} \leq 1\}.$$
The space $C^r(M)$ equipped with $\|\cdot\|_{r, \eps}$ is a Banach space. The induced topology is the standard $C^r$ topology.

Let $C^r(\M)$ be the Banach space (naturally identified with $C^r(M) \times C^r(M)$) of functions $$f : \RR^d \times \{0,1\} \rar \RR, (x,i) \mapsto f_i(x),$$ for which
 $f_i \in C^r(M),$ and equipped with the norm
 \beq
 \label{eq:normM}
 \|f\|_{r,\eps} = \|f_0\|_{r,\eps} + \|f_1\|_{r,\eps}.
 \eeq

\subsection{The induced chain and the $Q-\Delta$ decomposition}
From now on and until the end of Section \ref{sec:proof}, we assume that $r \in \NN$ is given and that the jump rates $\alpha_{01},\alpha_{10}$ satisfy condition (\ref{eq:maincond}).  We fix $\alpha, \eps > 0$ such that
\beq
\label{eq:maincondbis}
\alpha > \max(\alpha_{01},\alpha_{10}) \geq \min(\alpha_{01},\alpha_{10}) >  r |\Lambda_1| + \sum_{i = 1}^d |\Lambda_i| + r\eps.
\eeq
Let $U$ be the Markov transition matrix on $\{0,1\}$ given by
$$
U = \frac{1}{\alpha} \left(
                         \begin{array}{cc}
                           \alpha - \alpha_{01} & \alpha_{01} \\
                           \alpha_{10} & \alpha -\alpha_{10} \\
                         \end{array}
                       \right).
                       $$
It induces a Markov operator on $\M,$  still denoted $U$, defined as
$$U f(x,i) :=  (U f(x))_i = U_{i0} f(x,0) + U_{i1}f(x,1).$$
Let now, $K$ and $P$ be the Markov operators on $\M$  respectively defined by
$$K f(x,i) = \int_0^{\infty} f(\Phi_i^{t/\alpha}(x),i)e^{-t} dt,$$
and
\beq
\label{eq:defP}
P = KU.
\eeq
\brem
\label{rem:chainXI}
{\rm The kernel $P$ is the kernel of a discrete time Markov chain  $(X_n,I_n)_{n \geq 0}$ living on $\M$ where
$(I_n)$ is a Markov chain on $\{0,1\}$ independent of $(X_n)$ with transition matrix $U,$
 $$X_{n+1} = \Phi_{I_n}^{\tau_{n+1}}(X_n),$$
and  $(\tau_n)_{n \geq 0}$ is a sequence of i.i.d random variables having an exponential distribution with parameter $\alpha.$
 }
 \erem
The (unique) invariant distribution $\pi$ of $P$ and the (unique) invariant distribution $\Pi$ of the PDMP ((\ref{eq:pdmp}), (\ref{eq:jumps})) are linked by the following formulaes established in (\cite{BLMZ_2015}, Proposition 2.4 and Lemma 2.6).
$$
\Pi = \pi K, \, \pi =  \Pi U.
$$
In particular,  the densities of $\pi$ and $\Pi,$ satisfy
\beq
\label{eq:dpiPi}
\frac{d\pi_i}{dx} = \frac{d\Pi_0}{dx}U_{0i} + \frac{d\Pi_{1}}{dx}U_{1i}, i = 0,1.
 \eeq
 The $2 \times 2$ matrix $U$ being invertible, this shows that $\frac{d\pi}{dx}$ and $\frac{d\Pi}{dx}$ have the same regularity properties.
The proof of Theorem \ref{th:main} then amounts to proving  that $\pi$ has a $C^r$ density.

The strategy of the proof heavily relies on the recent paper  \cite{benaïm2024invariant}. The main idea is to show that, for $n$ sufficiently large, $P^n$ decomposes as the sum of a sub-markov kernel $Q$ producing regularity and a "small" sub-Markov kernel $\Delta$ preserving regularity.

Let ${\cal M}(\M)$   be the set of nonnegative bounded measures on $\M$ and ${\cal C}(\M) \subset {\cal M}(\M)$ the subset of measures $\mu$ that are absolutely continuous with respect to the Lebesgue measure $\m$ and whose  density $\rho$ lies in $C^r(\M).$ For such a $\mu$  we write
$\|\mu\|_{r,\eps}$ for $\displaystyle \|\rho\|_{r,\eps}$ by abuse of notation.

 A kernel $Q$ on $\M$ is called {\em subMarkov} (respectively {\em nondegenerate}) if  $Q(z,\M) \leq 1,$ (respectively  $Q(z, \M) \neq 0$) for all $z \in \M.$ It is called {\em Feller} if $Q f$ is continuous whenever $f$ is continuous on $\M.$

 The following lemma is the cornerstone of the proof. Combined with Theorem 2.8 in \cite{benaïm2024invariant}, it implies that $\pi,$ hence $\Pi,$ lies in ${\cal C}(M).$
\blem[Key Lemma]
\label{lem:key}
Assume that condition (\ref{eq:maincond}) holds.  Then, there exists $n \in \NN^*$   such that
$P^n = Q + \Delta,$ where $Q$ is a nondegenerate sub-Markov kernel and $\Delta$ a sub-Markov kernel on $\M$ satisfying the following properties:
\bdes
\iti ${\cal M}(\M) Q \subset {\cal C}(\M),$
\itii ${\cal C}(\M) \Delta \subset {\cal C}(\M)$
\itiii For all $\mu \in {\cal C}(\M), \sum_k \|\mu \Delta^k\|_{r,\eps} < \infty.$
\edes
\elem
The remainder of the section is devoted to the proof of this lemma.
\subsection{Proof of the key Lemma \ref{lem:key}}
Given $n \geq 1, \i = (i_1, \ldots, i_n) \in \{0,1\}^n,$ and $s = (s_1, \ldots, s_n) \in \RR^n,$  we set
\beq
\label{eq:defbigphi}
\mathbf{\Phi}^s_{\i} = \Phi_{i_n}^{s_n} \circ \ldots \circ \Phi_{i_1}^{s_1}.
\eeq
and
\beq
\label{eq:defUi}
U[\i] =   U_{i_1 i_2}  \ldots U_{i_{n-2} i_{n-1}} U_{i_{n-1} i_n}.
\eeq
If $i,j  \in \{0,1\}$  we use the notation $(i,\i), (i,\i,j),$ etc. to denote the sequences
$(i,i_1, \ldots, i_{n}) \in \{0,1\}^{n+1}, (i, i_1,\ldots,i_{n},j) \in \{0,1\}^{n+2},$ etc.

Given a  $C^{\infty}$ function  $h : (\RR_{+}^*)^{n} \mapsto [0,1]$ and $\i \in \{0,1\}^{n-1},$  we let  $P_{ \i, h}$ denote the sub-Markovian operator on $\M$
 defined by
\beq
\label{eq:defPih}
P_{ \i, h} f(x,i) = \sum_{j \in \{0,1\}} U[(i, \i, j)] \int_{\RR^{n}_+} f(\mathbf{\Phi}_{(i,\i)}^{s/\alpha}(x),j)   e^{-|s|}h(s) ds,
\eeq where
$|s|$ stands for  $s_1 + \ldots + s_{n}.$

 If $h \equiv 1,$  we write   $P_{ \i}$ for $P_{ \i, h.}$
Clearly we have that
$$P_{ \i} = P_{ \i, 1-h} + P_{ \i, h}$$
and
$$P_{ \i} f(x,i) = \EE_{x,i}(f(X_{n}, I_{n}) \Ind_{\left \{(I_1,\ldots, I_{n-1})  = \i \right \}}),$$ where $(X_n,I_n)$ is the discrete time Markov chain having $P$ as transition kernel (see Remark \ref{rem:chainXI}).
In particular
\beq
\label{eq:defPn}
P^{n}  = \sum_{\i \in \{0,1\}^{n-1}} P_{ \i}.
\eeq
For $n > 1$ and $\i \in \{0,1\}^n,$ let $$S(\i) = \mathsf{card}\{k \in \{1, \ldots, n-1\}\: : i_{k+1} \neq i_k\}$$ denote the number of "switches" in the sequence $\i.$
For all $n > d+1,$ we can then rewrite $P^n$ as $P^{n} = Q + \Delta,$ where  $\Delta = \Delta_1 + \Delta_2,$
\beq
\label{eq:Delta1}
\displaystyle \Delta_1 =  \sum_{\{\i \in \{0,1\}^{n-1} \: : \mathsf{S}(\i) < d\}} P_{ \i},
\eeq
\beq
\label{eq:Delta2}
\Delta_2 = \sum_{\{\i \in \{0,1\}^{n-1} \: : \mathsf{S}(\i) \geq d\}} P_{ \i, 1- h_{\i}},
\eeq
\beq
\label{eq:Q}
Q =  \sum_{\{\i \in \{0,1\}^{n-1} \: : \mathsf{S}(i) \geq d\}} P_{ \i, h_{\i}},
\eeq
 and  $h_{\i} : (\Rp^*)^n \rar  [0,1]$ (appearing in $\Delta_2$ and $Q$) is a smooth function to be chosen later.
 Our next goal is to show that $Q$ and $\Delta$ meet the requirement of Lemma \ref{lem:key}.
\subsubsection*{Tranfer operators}
The {\em transfer operator} of $\Phi^t_i$ is the operator   on $L^1(dx)$ defined as
\beq
\label{eq:transf} {\cal L}_{\Phi^t_i}(\rho)(x) = \rho(\Phi^{-t}_i(x)) \mathsf{det}(D\Phi^{-t}_i(x)) = \rho(\Phi^{-t}_i(x)) e^{- t \mathsf{Tr}(A)}
\eeq
where $$\mathsf{Tr}(A) = \Lambda_1 + \ldots + \Lambda_d$$ is the trace of $A.$

Note that if $\rho$ is the density (with respect to the Lebesgue measure) of some bounded measure $\mu,$  then ${\cal L}_{\Phi^t_i}(\rho)$ if the density of $\mu \circ (\Phi_i^t)^{-1},$ the image measure of $\mu$ by $\Phi_i^t.$

Note also that for all $t \geq 0$
$${\cal L}_{\Phi^t_i}(C^r(M))  \subset C^r(M).$$

The following lemma follows directly from the chain rule and the definition of $\|\cdot\|_{r, \eps}$ on $C^r(M).$ Here $\|{\cal L}_{\Phi_i^t}\|_{r,\eps}$ stands for the operator norm of ${\cal L}_{\Phi_i^t}$ on $(C^r(M),\, \|\cdot\|_{r,\eps}).$
\blem For all $t \geq 0,$ and $i \in \{0,1\},$
$$\|{\cal L}_{\Phi_i^t}\|_{r,\eps} \leq \exp{ \left[-t \left( \mathsf{Tr} A - r (|\Lambda_1| + \eps) \right)\right]} =
\exp{\left[t \left(r|\Lambda_1| + \sum_{i = 1}^d |\Lambda_i| + r \eps \right)\right]}.$$
\elem
For all $\i \in \{0,1\}^n$ and $s \in \RR^n,$ the   transfer operator of $\mathbf{\Phi}^s_{\i}$ (see equation (\ref{eq:defbigphi})) on $L^1(dx)$  is naturally defined as
$${\cal L}_{\mathbf{\Phi}^s_{\i}} = {\cal L}_{\Phi_{i_n}^{s_n}} \circ \ldots \circ  {\cal  L}_{\Phi_{i_1}^{s_1}}$$
where ${\cal  L}_{\Phi_{i}^{t}}$ is defined by (\ref{eq:transf}).

Associated to $P_{ \i, h}$ (see equation (\ref{eq:defPih}))  is the transfer operator  on $L^1(\m)$  given by
\beq
\label{eq:transP}
{\cal P}_{\i,h}(\rho)(x,j) = U[0,\i,j]{\cal L}_{(0,\i),h}(\rho_0) +  U[1,\i,j]{\cal L}_{(1,\i),h}(\rho_1)
\eeq
where
\beq
\label{eq:transL}
{\cal L}_{(i,\i),h}(\rho_i) := \int_{\Rp^n} {\cal L}_{\mathbf{\Phi}_{(i,\i)}^{s/\alpha}}(\rho_i) e^{-|s|} h(s) ds
\eeq
Here also, when   $h \equiv 1,$  we write ${\cal P}_{\i}$ for ${\cal P}_{\i,h}.$

It is readily seen that if $\mu \in  {\cal M}(\M)$ has density $\rho \in L^1(\m)$ then $\mu P_{\i,h}$ has density ${\cal P}_{\i,h}(\rho).$
Furthermore,
\blem
\label{lem:boundPih} The operator  ${\cal P}_{\i,h}$ is bounded on $C^r(\M)$ and for every set $I \subset \{0,1\}^{n-1},$ letting
$U[(i,I)] =  \sum_{\i \in I} U[(i,\i)],$ one has
$$\displaystyle \|\sum_{\i \in I} {\cal P}_{\i, h}\|_{r,\eps} \leq \max(U[(0,I)], U([(1,I)]) C(h)$$
where
$$C(h)=\int_{(\Rp^*)^n} \exp |s| \left[ \frac{r|\Lambda_1| + \sum_{i = 1}^d |\Lambda_i| + r \eps }{\alpha} - 1 \right] h(s) ds$$
$$  \leq \left[ 1- \frac{r|\Lambda_1| + \sum_{i = 1}^d |\Lambda_i| + r \eps }{\alpha}\right]^{-n}.$$
\elem
\prf Set $C = C(h).$ By the preceding lemma and the definition (\ref{eq:transL}) one has
$$\|{\cal L}_{(i,\i),h}(\rho_i)\|_{r,\eps} \leq C(h) \|\rho_i\|_{r,\eps}.$$ Thus, by the definition of the norm $\| \cdot \|_{r,\eps}$ on $\M$ it comes
that $$\|\sum_{\i \in I} {\cal P}_{\i, h}(\rho)\|_{r,\eps} \leq C \sum_{\i  \in I} (U[(0,\i,0)]  \|\rho_0\|_{r,\eps} + U[(1,\i,0)]  \|\rho_1\|_{r,\eps}) $$
$$+  C \sum_{\i  \in I} (U[(0,\i,1)]  \|\rho_0\|_{r,\eps} + U[(1,\i,1)]  \|\rho_1\|_{r,\eps}) $$
$$\leq \max(U[(0,I)], U([(1,I)]) C (\|\rho_0\|_{r,\eps} +  \|\rho_1\|_{r,\eps}) = \max(U[(0,I)], U([(1,I)]) C \|\rho\|_{r,\eps}.$$
This proves the first inequality. The last one is immediate because $0 \leq h \leq 1.$
\qed
\subsubsection*{The operator $\Delta_1$}
A consequence of the preceding  lemma is that, under the assumption of the key Lemma, the norm of the transfer operator associated to $\Delta_1$ (equation (\ref{eq:Delta1})) can be made arbitrary  small by choosing $n$ sufficiently large.
\blem
\label{lem:boundPi} For all $\eta > 0$ there exists $n$ sufficiently large  such that
$$\|\sum_{\{\i \in \{0,1\}^{n-1} \: : \mathsf{S}(\i) < d\}} {\cal P}_{\i}\|_{r,\eps} < \eta.$$
\elem
\prf Assume for simplicity that $d$ is even (the proof is the same with $d$ odd). Let $T^1, T^2, \ldots$ be independent random variables taking values in $\NN^*$ such that, for $i$ odd (respectively even)
$T^i$ has a geometrical distribution with parameter $U_{01}$ (respectively $U_{10}$).
Then,
$$\sum_{\{\i \in \{0,1\}^{n-1} \mathsf{S}(\i) < d\}} U(0,\i) = \Pr(T^1+ \ldots + T^d > n-1).$$
Now, by Markov inequality, for all
\beq
\label{eq:condz}
 1 < z < \frac{1}{\max(U_{00},U_{11})} = \min ((1-\frac{\alpha_{01}}{\alpha})^{-1},(1-\frac{\alpha_{10}}{\alpha})^{-1}),
 \eeq
 $$ \Pr(T^1+ \ldots + T^d > n-1) \leq \frac{\E(z^{T^1+ \ldots +T^d})}{z^{n-1}} =  \frac{1}{z^{n-1}}[\frac{z(1-U_{00})}{(1-z U_{00})}]^{d/2} [\frac{z(1-U_{11})}{(1-z U_{11})}]^{d/2}.$$
 Our condition  \ref{eq:maincondbis} permits to choose  $z$ satisfying (\ref{eq:condz}) such that
 $$z \left[ 1- \frac{r|\Lambda_1| + \sum_{i = 1}^d |\Lambda_i| + r \eps }{\alpha}\right] > 1.$$ Combined with Lemma \ref{lem:boundPih}, this proves the result.
 \qed
\subsubsection*{The operators $Q, \Delta_2,$ and the end of the proof}
  For $\i \in \{0,1\}^n,$ we use the notation $\mathbf{\Phi}_{\i}(x)$ to denote the mapping $s \in \RR^n \mapsto \mathbf{\Phi}^{s}_{\i}(x) \in \RR^d$  where $\mathbf{\Phi}^{s}_{\i}(x)$ is defined by (\ref{eq:defbigphi}).
\blem
\label{lem:subm1} Let $\mathbf{j} \in \{0,1\}^{d+1}$ be one of the two  sequences defined by $j_{k+1} = 1-j_k$ for $k = 1, \ldots, d.$ That is
$\mathbf{j}= (0,1,0,\ldots)$ or $\mathbf{j} = (1,0,1, \ldots).$ Then, there exists an open set
 $O \subset (\Rp^*)^{d+1}$ having full Lebesgue measure such that {\bf for all} $x \in M,$ and $s \in O, \mathbf{\Phi}_{\mathbf{j}}(x)$ is  submersive at $s.$
\elem
\prf

Assume for the sake of the proof that $j_1 = 0.$ The other case $j_1 = 1$ is treated similarly.
Set $\mathbf{\Phi}^s = \mathbf{\Phi}_{\mathbf{j}}^s(0).$
Let $$O' = \{(s_2, \ldots, s_{d+1}) \in (\Rp^*)^{d} \: : \mathsf{det}\left(\partial_{s_2}\mathbf{\Phi}^s, \ldots, \partial_{s_{d+1}}\mathbf{\Phi}^s \right) \neq 0\}.$$
Our first goal is to show that $O'$ has full measure in $(\Rp^*)^{d}.$

To shorten notation, write $y_k = \Phi_{j_k}^{s_k} \circ \ldots \circ \Phi_{j_1}^{s_1}(0).$
 Using the chain rule, it comes that, for $k = 1, \ldots, d+1,$
$$ \partial_{s_k} \mathbf{\Phi}^s = D(\Phi_{j_{d+1}}^{t_{d+1}} \circ \dots \circ \Phi_{j_{k+1}}^{t_{k+1}})(y_k) F_{j_k}(y_k) $$
$$ = e^{\left( (\sum_{i = k+1}^{d+1} s_i)  A\right)} F_{j_k}(y_k).$$
Thus, $\partial_{s_1} \mathbf{\Phi}^s = 0,$ and for all $k = 1, \ldots, d+1,$
$$ \partial_{s_k} \mathbf{\Phi}^s -  \partial_{s_{k-1}} \mathbf{\Phi}^s
=  e^{\left( (\sum_{i = k+1}^{d+1} s_i)  A\right)} [F_{j_k}(y_k) - e^{s_k A} F_{j_{k-1}}(y_{k-1})]$$
$$=  \pm A e^{ (\sum_{i = k+1}^{d+1} s_i)  A)} p.$$ Here, the last equality follows from the simple observation that
$$F_{i}(\Phi_{1-i}^s(y)) - e^{s A} F_{1-i}(y) = (-1)^i Ap$$ for all $y \in \RR^d, s \in \RR$ and $i \in \{0,1\}.$
 Rewriting
 $$\mathsf{det}\left(\partial_{s_2}\mathbf{\Phi}^s, \ldots, \partial_{s_{d+1}}\mathbf{\Phi}^s \right) = \mathsf{det}\left(\partial_{s_2}\mathbf{\Phi}^s, \partial_{s_3}\mathbf{\Phi}^s - \partial_{s_2}\mathbf{\Phi}^s, \ldots,
\partial_{s_{d+1}}\mathbf{\Phi}^s - \mathbf{\Phi}(s_d)\right),$$ it follows that
$O'= (\Rp^*)^d \setminus (f \circ T)^{-1}(0),$ where  $$T(s_2, ..., s_{d+1}) = (s_2, s_2 + s_3, \ldots, s_2 + s_3 + \ldots + s_{d+1})$$  and
$$f(t_2, ..., t_{d+1}) = \mathsf{det}\left(e^{t_dA}p, ..., e^{t_2A}p\right).$$
The map $f$ being real analytical (by analyticity of $\mathsf{det}$ and $\mathsf{exp}$), $f^{-1}(0)$ has  measure zero unless $f$ is identically zero.  A proof of this fact is given for convenience   in the appendix (see Lemma  \ref{lem:analytic}). To see that $f$ is nonzero, notice that
$\partial_{t_2} \partial^2_{t_3^2} \ldots \partial^d_{t_{d+1}^d}f(0) = \mathsf{det}(A^d p, \ldots, A^2 p, Ap) \neq 0$
by Assumption \ref{mainhyp}. Now, since $T$ preserves the Lebesgue measure, $(f \circ T)^{-1}(0)$ has measure zero. This concludes the proof that $O'$ has full measure.

To conclude the proof of the lemma, observe that for all $x \in \RR^d,$ and $s \in (\Rp^*)^{d+1},$
 $\mathbf{\Phi}^s(x) = \mathbf{\Phi}^s + e^{(\sum_{i = 1}^{s_{d+1}} s_i A)} x.$
Therefore,
$$\partial_{s_k} \mathbf{\Phi}^s(x)- \partial_{s_1} \mathbf{\Phi}^s(x) = \partial_{s_k} \mathbf{\Phi}^s$$ for all $k = 2, \ldots d+1.$
This shows that for all $s \in O := \Rp^* \times O'$ and for all $x \in \RR^d,$ $\mathbf{\Phi}(x)$ is submersive at $s.$
\qed
The following lemma follows directly from the preceding one.
\blem
\label{lem:subm2} Let $n \geq d+1$ and $\i \in \{0,1\}^n$ be such that $S(\i) \geq d.$
 \bdes
 \iti There exists an open set $O_{\i} \subset (\Rp^*)^n$ having full Lebesgue measure such that for all $s \in O_{\i}$ and $x \in M,$ $\mathbf{\Phi}_{\i}(x)$ is submersive at $s;$
 \itii For all $\eta > 0$ there exists a $C^{\infty}$ function $h_{\i} : (\Rp^*)^n \rar [0,1],$ with support in $O_{\i},$ such that
$$\|{\cal P}_{\i,1-h_{\i}}\|_{r,\eps} \leq \eta;$$ and
$${\cal M}(\M) P_{\i, h_{\i}} \subset {\cal C}(\M).$$
\edes
\elem
\prf
The condition  $S(\i) \geq d,$ implies that
$$\mathbf{\Phi}_{\i}^{(s_1, \ldots, s_n)}(x) = \Phi_{i_n}^{s_n} \circ \dots \circ \Phi_{i_{k + 1}}^{s_{k + 1}} (\mathbf{\Phi}_{\j}^{T(s_1, \ldots, s_{k})}(x))$$ where $\j \in \{0,1\}^{d+1}$ is like in Lemma \ref{lem:subm1}, $d+1 \leq k \leq n,$ and $T : \RR^k \rar \RR^{d+1}$ writes
$$T(s_1, \ldots, s_{k}) = (\sum_{0 < i \leq k_1} s_i, \sum_{k_1 < i  \leq k_2} s_i,   \ldots, \sum_{k_d < i  \leq k} s_i).$$
with $1  \leq k_1 < k_2 < \ldots < k := k_{d+1} \leq n.$
By Lemma \ref{lem:subm1}, and surjectivity of $T,$ $s \mapsto \mathbf{\Phi}_{\j}^{T(s)}(x)$ is a submersion on $T^{-1}(O)$ where $O$ is like in Lemma \ref{lem:subm1}. Since
 $x \mapsto \Phi_{i_n}^{s_n} \circ \dots \circ \Phi_{i_{k + 1}}^{s_{k + 1}}(x)$ is a diffeomorphism,
 $\mathbf{\Phi}_{\i,x}$ is then a submersion at every $s  \in O_{\i} = T^{-1}(O) \times  (\Rp^*)^{n-(k + 1)}.$ Furthermore,  $O_{\i}$ has full measure in $\RR^n$ because $O$ has full measure and $T$ is a linear surjection. This proves $(i).$

 We now prove $(ii).$ Let $\nu$ be the Borel measure on $\RR^n$ defined by   $$\nu(ds) = \exp |s| \left[ \frac{r|\Lambda_1| + \sum_{i = 1}^d |\Lambda_i| + r \eps }{\alpha} - 1 \right] \Ind_{ (\Rp^*)^n}(s).$$ Then, $\nu$ is finite, hence regular. Therefore,  there exists a compact set $K_{\i} \subset O_{\i}$ such  that $\nu(O_{\i} \setminus K_{\i}) \leq \eta.$ Let $h_{\i}$ be a $C^{\infty}$ function with $0 \leq h_{\i} \leq 1$ having support in $O_{\i}$ and such $h|_{K_{\i}} = 1.$ Then the first assertion follows from Lemma \ref{lem:boundPih} and the second from Proposition 3.2 in \cite{benaïm2024invariant}, because by $(i)$   $\mathbf{\Phi}_{\i}(x)$ is submersive at every point $s$ in the support of $h_{\i}.$
\qed
We can now conclude the proof of Lemma \ref{lem:key}. Fix $\eta < \frac{1}{2}.$ First, choose $n$ such that the conclusion  of Lemma \ref{lem:boundPi} holds.
Then, relying on the preceding lemma  \ref{lem:subm2},   choose  for each $\i \in \{0,1\}^{n-1},$ a smooth function $0 \leq h_{\i} \leq 1$ such that
$\|{\cal P}_{\i,h_{\i}}\|_{r,\eps} \leq \frac{\eta}{2^n}.$ Next define $Q$ and $\Delta = \Delta_1 + \Delta_2$  by formulaes  (\ref{eq:Q}), (\ref{eq:Delta1}) and (\ref{eq:Delta2}).

\section{Appendix}
The Lebesgue measure on $\RR^m$ is denoted $\mathsf{Leb}_m.$
\blem
\label{lem:analytic}
Let $m \geq 1$ and  let $f : I_1 \times \ldots \times I_m \rar \RR$ be a real analytic function, where $I_i \subset \RR, i = 1, \ldots, m$ are  nonempty open intervals. Then  either $f = 0$ (i.e $f$ is the zero function) or $\mathsf{Leb}_m(f^{-1}(0)) = 0.$
\elem
\prf We proceed by induction on $m.$ If $m = 1,$ and $f \neq 0,$  $f^{-1}(0)$ is countable by analyticity. Hence $\mathsf{Leb}_1(f^{-1}(0)) = 0.$ Suppose that the property is true for $m-1$ with $m \geq 2.$ Set $J = I_2 \times \ldots \times I_m.$ Then, by Fubini-Tonelli,
$$\mathsf{Leb}_m(f^{-1}(0)) = \int_{I_1} \left (\int_J \Ind_{f^{-1}(0)}(x,y) dy \right ) dx = \int_{I_1} \mathsf{Leb}_{m-1}(f_x^{-1}(0)) dx$$
where $f_x : J \rar \RR$ is the map defined by $f_x(y) = f(x,y).$

Let $A = \{x \in I_1: \: f_x = 0\}.$ Then, by continuity, $$A = \bigcap_{y \in J \cap \mathbb{Q}^{m-1}} \{x \in I_1: \: f(x,y) = 0\}. $$ This shows that $A$ is measurable. Hence, by the induction hypothesis,
$$\mathsf{Leb}_m(f^{-1}(0)) = \int_A \mathsf{Leb}_{m-1}(f_x^{-1}(0)) dx.$$
Now, if $f \neq 0$ there exists some $x^* \in I_1$ and $y^* \in J \cap \mathbb{Q}^{m-1}$ such that $f(x^*,y^*) \neq 0.$ It comes that
$$\int_A \mathsf{Leb}_{m-1}(f_x^{-1}(0)) dx \leq \int_{\{x \in I_1: \: f(x,y^*) = 0\}} \mathsf{Leb}_{m-1}(f_x^{-1}(0)) dx $$
This later integral is zero because, by the base case, $$\mathsf{Leb}_1(\{x \in I_1: \: f(x,y^*) = 0\}) = 0.$$
\qed

 {\textbf{Acknowledgement:}}  The work of MB,  was funded by the grant 200020-219913 from the Swiss National Foundation. We thank Yuri Bakhtin, Tobias Hurth, Claude Lobry, Edouard Strickler and Oliver Tough for stimulating discussions.

\bibliographystyle{amsplain}

\bibliography{regularlineardensity.bib}

\providecommand{\bysame}{\leavevmode\hbox to3em{\hrulefill}\thinspace}
\providecommand{\MR}{\relax\ifhmode\unskip\space\fi MR }
\providecommand{\MRhref}[2]{%
  \href{http://www.ams.org/mathscinet-getitem?mr=#1}{#2}
}
\providecommand{\href}[2]{#2}
\begin{thebibliography}{10}

\bibitem{ArnoldKliemann}
L.~Arnold and W.~Kliemann, \emph{On unique ergodicity for degenerate
  diffusions}, Stochastics \textbf{21} (1987), no.~1, 41--61. \MR{899954}

\bibitem{bakhtin&hurt&matt}
Y.~Bakhtin, T.~Hurth, and J.C. Mattingly, \emph{Regularity of invariant
  densities for 1d-systems with random switching}, Nonlinearity \textbf{28}
  (2015), 3755--3787.

\bibitem{bakhtin&hurt}
Yuri. Bakhtin and Tobias. Hurth, \emph{Invariant densities for dynamical
  systems with random switching}, Nonlinearity (2012), no.~10, 2937--2952.

\bibitem{BHLM18}
Yuri Bakhtin, Tobias Hurth, Sean~D. Lawley, and Jonathan~C. Mattingly,
  \emph{Smooth invariant densities for random switching on the torus},
  Nonlinearity \textbf{31} (2018), no.~4, 1331--1350. \MR{3816636}

\bibitem{BHLM21}
\bysame, \emph{Singularities of invariant densities for random switching
  between two linear {ODE}s in 2{D}}, SIAM J. Appl. Dyn. Syst. \textbf{20}
  (2021), no.~4, 1917--1958. \MR{4322095}

\bibitem{BHS18}
M.~Bena\"im, T.~Hurth, and E.~Strickler, \emph{A user-friendly condition for
  exponential ergodicity in randomly switched environments}, Electronic
  Communications in Probability \textbf{23} (2018), no.~44, 1--12.

\bibitem{BLMZ_2015}
M.~Bena\"{\i}m, S.~Le~Borgne, F.~Malrieu, and P.-A. Zitt, \emph{Qualitative
  properties of certain piecewise deterministic {M}arkov processes}, Ann. Inst.
  Henri Poincar\'{e} Probab. Stat. \textbf{51} (2015), no.~3, 1040--1075.

\bibitem{BLMZ_2019}
\bysame, \emph{Erratum: Qualitative properties of certain piecewise
  deterministic {M}arkov processes}, Ann. Inst. Henri Poincar\'{e} Probab.
  Stat. \textbf{55} (2019), no.~4.

\bibitem{BenaimHurth}
Michel Bena\"{\i}m and Tobias Hurth, \emph{Markov chains on metric spaces, a
  short course}, Universitext, vol.~99, Springer, Cham, 2022.

\bibitem{BenaimLobry}
Michel Bena\"{\i}m and Claude Lobry, \emph{Lotka-{V}olterra with randomly
  fluctuating environments or ``{H}ow switching between beneficial environments
  can make survival harder''}, Ann. Appl. Probab. \textbf{26} (2016), no.~6,
  3754--3785. \MR{3582817}

\bibitem{BCL17}
M~Benaïm, F~Colonius, and R~Lettau, \emph{Supports of invariant measures for
  piecewise deterministic {M}arkov processes}, Nonlinearity \textbf{30} (2017),
  no.~9, 3400.

\bibitem{BT23}
Michel Benaïm and Oliver Tough, \emph{Regularity of the stationary density for
  systems with fast random switching}, 2023, arXiv 2212.03632.

\bibitem{benaïm2024invariant}
\bysame, \emph{On invariant distributions of feller markov chains with
  applications to dynamical systems with random switching}, 2024, arXiv
  2310.17543.

\bibitem{MR3983339}
Joris Bierkens, Gareth~O. Roberts, and Pierre-Andr\'e Zitt, \emph{Ergodicity of
  the zigzag process}, Ann. Appl. Probab. \textbf{29} (2019), no.~4,
  2266--2301. \MR{3983339}

\bibitem{Bony}
Jean-Michel Bony, \emph{Principe du maximum, in\'{e}galite de {H}arnack et
  unicit\'{e} du probl\`eme de {C}auchy pour les op\'{e}rateurs elliptiques
  d\'{e}g\'{e}n\'{e}r\'{e}s}, Ann. Inst. Fourier (Grenoble) \textbf{19} (1969),
  no.~fasc. 1, 277--304 xii. \MR{262881}

\bibitem{Bressloff17}
Paul~C. Bressloff, \emph{Stochastic switching in biology: from genotype to
  phenotype}, J. Phys. A \textbf{50} (2017), no.~13, 133001, 136. \MR{3623058}

\bibitem{Cloez15}
Bertrand Cloez and Martin Hairer, \emph{Exponential ergodicity for {M}arkov
  processes with random switching}, Bernoulli \textbf{21} (2015), no.~1,
  505--536. \MR{3322329}

\bibitem{MR3582809}
Manon Costa, \emph{A piecewise deterministic model for a prey-predator
  community}, Ann. Appl. Probab. \textbf{26} (2016), no.~6, 3491--3530.
  \MR{3582809}

\bibitem{Davis84}
M.~H.~A. Davis, \emph{Piecewise-deterministic {M}arkov processes: a general
  class of nondiffusion stochastic models}, J. Roy. Statist. Soc. Ser. B
  \textbf{46} (1984), no.~3, 353--388, With discussion. \MR{790622}

\bibitem{Lawley15}
Sean~D. Lawley, Jonathan~C. Mattingly, and Michael~C. Reed, \emph{Stochastic
  switching in infinite dimensions with applications to random parabolic
  {PDE}}, SIAM J. Math. Anal. \textbf{47} (2015), no.~4, 3035--3063.
  \MR{3379020}

\bibitem{li17}
Dan Li, Shengqiang Liu, and Jing'an Cui, \emph{Threshold dynamics and
  ergodicity of an {SIRS} epidemic model with markovian switching}, Journal of
  Differential Equations \textbf{263} (2017), no.~12, 8873--8915.

\bibitem{Malrieu}
Florent Malrieu, \emph{Some simple but challenging {Markov} processes}, Annales
  de la Facult\'e des sciences de Toulouse : Math\'ematiques \textbf{Ser. 6,
  24} (2015), no.~4, 857--883 (en).

\bibitem{MR3695481}
Florent Malrieu and Pierre-Andr\'e Zitt, \emph{On the persistence regime for
  {L}otka-{V}olterra in randomly fluctuating environments}, ALEA Lat. Am. J.
  Probab. Math. Stat. \textbf{14} (2017), no.~2, 733--749. \MR{3695481}

\bibitem{Monemvassitis_2023}
Athina Monemvassitis, Arnaud Guillin, and Manon Michel, \emph{Pdmp
  characterisation of event-chain monte carlo algorithms for particle systems},
  Journal of Statistical Physics \textbf{190} (2023), no.~3.

\bibitem{Yin}
G.~G. Yin and C.~Zhu, \emph{Hybrid switching diffusions}, Stochastic Modelling
  and Applied Probability, vol.~63, Springer, New York, 2010, Properties and
  applications. \MR{2559912 (2010i:60226)}

\end{thebibliography}
\end{document}